%&amstex
\documentstyle {amsppt}
\magnification=\magstephalf

\def\sp{\  } \def\slash{ /} \def\a{ \alpha} \def\b{ \beta} \def\g{ \gamma}

\def\d{ \delta}    \def\z{
\zeta}
     \def\k{
\kappa} \def\l{ \lambda} \def\L{ \Lambda}         

\def\si{ \sigma} \def\Si{ \Sigma}

\def\o{ \omega}  \def\bq{ \Bbb Q}  \def\bz{ \Bbb
Z} \def\bc{
\Bbb C}  \def\8{_{ \infty}   } 
  \def\ps{ $p_1$-structure\sp} 
 \def\cy{ \bz[\l_p]}

\hfuzz 27pt

\topmatter \title On the Witten-Reshetikhin-Turaev representations of mapping
class groups
\endtitle
\rightheadtext{WRT representations}    \author Patrick M. Gilmer \endauthor
 \affil Louisiana State University \endaffil
 \address Department of Mathematics,  Baton Rouge, LA 70803 U.S.A  \endaddress
  \email gilmer\@ math.lsu.edu \endemail

\abstract

We consider a central extension of the mapping class group of a surface with a
collection of
framed colored points. The Witten-Reshetikhin-Turaev TQFTs associated to
$SU(2)$ and $SO(3)$
induce  linear representations of this group. We show that the denominators of
matrices
which describe these representations over a cyclotomic field can be restricted
in many
cases.  In this way, we give a proof of the known result that if the surface is
a torus with
no colored points, the representations have finite image. \endabstract \thanks
This research
was partially supported by a grant from the Louisiana Education Quality Support Fund
\endthanks
 \subjclass 57M99
\endsubjclass
  \endtopmatter \document

Recall that an object in a cobordism category of dimension 2+1 is a  closed
oriented
surfaces $\Si$ perhaps with some specified further structure. A morphism $M$
from $\Si$ to
$\Si'$ is (loosely speaking) a compact oriented 3-manifold perhaps with some
specified
further structure, called a cobordism, whose  boundary is  the disjoint union
of $-\Si$ and
$\Si'.$ A morphism $M'$ from $\Si'$ to $\Si''$ is composed with a morphism from
$\Si$ to
$\Si'$ by gluing along $\Si',$ inducing any required extra structure from the
structures on $M$ and $M'.$ Also the extra structure on a 3-manifold must
induce the extra structure on the boundary. A TQFT in dimension 2+1 is
then a functor from such a cobordism category to the category of modules over
some ring $R.$
There are further axioms that are generally required \cite{A,BHMV,Q}. One
usually denotes
the module associate to $\Si$ by $V(\Si),$ and denotes the homomorphism
associated to
cobordism $M$ by $Z(M).$   A TQFT yields a representation of (an extension) of
the mapping
class group of a surface. An extension is needed if there is some choice in the
extra data
which may be placed on a mapping cylinder.

We will study a version of the Witten-Reshetikhin-Turaev TQFTs \cite{W,RT}
associated to
$SU(2)$ and $SO(3) $ constructed by \cite{BHMV}.  In particular, we will use
the notation
where $V_p$ for $p=2r$ is  a TQFT associated to $SU(2).$  Also $V_p$ for $p$
odd is a TQFT
associated to $SO(3).$  We will assume that $p\ge3.$   We will  use a variant
of
the \cite{BHMV} approach  obtained by adapting an idea of Walker \cite{Wa}. We
replace, in the
definition of the cobordism category,  a \ps on a surface $\Si$ with  a
Lagrangian subspace
of $H_1(\Si,\Bbb Z),$ and a \ps on a 3-manifold $M$ by an integer and  a
Lagrangian
subspace for
the boundary of $M.$ If one does this, one may work over the ring $r_p=\bz
[A_p,\frac 1
p,u_p],$ where $ A_p$ is a primitive $2p$th root of unity, and $u_p^2=
A_p^{-6-\frac{p(p+1)}{2}}.$ So $u_p$ plays the role of $\k_p^3.$ This makes no
real difference
to the arguments below but simplifies some expressions. We prefer not to take
further roots
of unity than necessary. Let $\l_p$ be a root of unity such that $\bz[A_p,u_p]=
\cy.$ Thus $r_p=\bz[\l_p,\frac 1 p].$

An object in the cobordism category $\Cal C_p$ that we consider will be a
closed oriented
surface  (possibly empty) $\Si$  together with a  choice of Lagrangian subspace
$\ell \subset H_1(\Si,\Bbb Z),$ and
also with a (possibly empty) collection of banded points colored with 
p-colors.  A banded point is an oriented arc through the point. Here a p-color
is
an integer
from zero to $p\slash 2-2$ if $p$ is even. If $p$ is odd, a p-color is an even
nonnegative
integer less than or equal to $p-3.$  A morphism from $\Si$ to $\Si'$ is a
compact oriented
3-manifold $M$  whose boundary comes equipped with an identification with $-\Si
\coprod
\Si'$ where $\Si$ and $\Si'$ are possibly empty objects. Moreover the three
manifold  includes a (possibly empty)  colored fat graph. Here a fat graph is
an
oriented
surface which deformation retracts to a  trivalent graph which meets the
boundary along arcs
corresponding to the banded points.  A colored fat graph is a fat
graph whose core is given a $p$-admissible coloring \cite{BHMV}  with matching
colors
at the banded points. The 3-manifold is also
equipped with an integer $n(M)$ called the weight.  One defines $n( M' \circ
M)$ to be $n(M')+ n(M) +\si(
N(M)_*(\ell),\ell', N(M')^*(\ell'')).$ Here $\si$ denotes Wall's non-additivity
function
\cite{Wall,Wa} or (minus) the Maslov index, and  $N(M)_*: \L(H_1(\Si))
\rightarrow
\L(H_1(\Si')$  and $N(M')^*: \L(H_1(\Si'')) \rightarrow \L(H_1(\Si')$ denote
Lagrangian
actions induced by $M$ and $M'$  \cite{T}.
 Also $\ell,$ $\ell',$ and  $\ell''$ denote the
Lagrangian subspaces with which $\Si,$$\Si',$ and $\Si''$ are equipped.

Suppose $M$ is a closed morphism (i.e a morphism from $\emptyset$ to
$\emptyset$) in $\Cal C_p.$ $M$ may
be described by framed surgery along a framed link $L$ in $S^3$  which misses
the fat graph
$G$. Let $L(\o_p)$ denote the linear combination of framed links obtained by
replacing each
component by the skein element $\o_p$ (we use the notation of \cite{BHMV}). Let
$e(G)$ be the
linear combination of framed links obtained by expanding the graph on the
surface. Then
define $<M>_p\in r_p $ to be  $
u_p^{n(M)-\text{Sign}(L)} \eta_p $
times the Kauffman bracket polynomial of $e(G) \coprod  L(\o_p)$ evaluated at
$A=A_p.$
Here $\text{Sign}(L)$ denotes the signature of the linking matrix associated to
the framed link $L$, and $\eta_p \in r_p$ is the element defined at the beginging of \S 1.  The
universal construction in  \cite{BHMV} then yields a TQFT on $\Cal C_p$ over
the ring $r_p.$  Thus $V_p(\emptyset)= r_p$, and
$Z_p(M): V_p(\emptyset) \rightarrow V_p(\emptyset)$ is multiplication by
$<M>_p.$

Consider a closed connected oriented surface $F$  with a collection of banded
colored points. The mapping class group  of $F$ consists of isotopy classes
of
diffeomorphisms $f$ of
$F$ to itself preserving the orientation, the coloring, and the framing on
the points. Thus points may  be permuted by $f$ if they happen to be colored
the same. Now equip $F$ with a Lagrangian subspace $\ell.$ Thus $(F,\ell)$ may
be considered as an object $\Si$ of $\Cal C_p,$ as long as the colors on $F$
are p-colors.
Consider $\Cal M(\Si),$ the central extension of the mapping class group of $F$
 whose
elements are pairs $(f,n),$ where $n$ is an integer.   If $(f_1,n_1):(F,\ell)
\rightarrow (F,\ell),$ and
$(f_2,n_2):(F,\ell) \rightarrow (F,\ell ),$ then $(f_2,n_2) \circ
(f_1,n_1)$ is defined by
$\left(f_2 \circ f_1,n_1 + n_2 +\si(f_2 \circ f_1(\ell),
f_2(\ell),\ell)\right).$

As an object of $\Cal C_p$, $\Si$   is assigned a free $r_p$-module
$V_p(\Si),$ which comes equipped with a
unimodular
Hermitian form $<\ ,\  >_\Si.$ The mapping cylinder of $f$ weighted by an
integer $n,$ ${C(f,n)},$ is a morphism from $\Si$ to itself in $\Cal C_p$. Thus
${C(f,n)}$ induces an endomorphism
$Z_p({C(f,n)})$ of $V_p(\Si).$ As ${C(f,n)}$ is a mapping cylinder,
$Z_p({C(f,n)})$ is an
isometry of $<\ ,\  >_\Si.$  This would not be true for an arbitrary morphism
in the
cobordism category. In this way, we obtain a representation $\rho_p$ of $\Cal
M(\Si)$ into
the group of isometries of $V_p(\Si).$

Let $H$ be a handlebody weighted zero with boundary $\Si$ such that the kernel
of the map on first
homology induced by the inclusion is the Lagrangian subspace with which $\Si$
is equipped.
    Let $\Cal G$ be  a fat graph in $H$ with boundary the framed
points, such that $H$ deformation retracts to $\Cal G.$   The $p$-admissible
colorings
$\b_i$ of $\Cal G$ which extend the given coloring on the boundary describe a
basis $\Cal B$
for $V_p(\Si),$ \cite{BHMV,4.11}. This basis is orthogonal with respect to the
unimodular
form defined on $V_p(\Si).$ In the case of a torus with no colored points ,
$\Cal B$ is
orthonormal. We have two overlapping results.

\proclaim{Theorem 1} Let $\Si$ have genus $g,$ and $p$ be either an
odd prime,
or  twice an odd prime. With respect $\Cal B,$ the matrices for $\rho_p$
have entries
lieing in $\left({\frac 1 p}\right)^{[\frac{g+1}{2}]}\cy.$
 \endproclaim

Here $[\ \ ]$ denotes the greatest integer function. In the lemmas below, we
specify a
somewhat stronger and more general restrictions.

\proclaim{Theorem 2} If $\Si$ is a torus without any colored points and $p$ is
even, the matrices
for $\rho_p$
have entries lieing in $\frac 1 {p} \cy.$
\endproclaim

 We find these theorems intriguing
as one would
expect in general that denominators should increase when multiplying matrices.

We obtain the following Corollary which M. Kontsevich informs us is known.
Kontsevich
observed it in 1988. It is also known in the conformal field theory community.
We thank M.
Kontsevich and G. Masbaum for pointing out an error in a previous version of
this paper.
We also thank the referee for useful comments.

\proclaim{Corollary} If $\Si$ is a torus without any colored points, then
$\rho_p$ has
finite image in  the isometries of $V_p(\Si).$ \endproclaim

\demo{Proof} Using \cite{BHMV,1.5}, one may deduce the result for $p$ odd from
the result
for $p$ even. So we assume $p$ is even.  $\cy$ maps to a discrete subgroup of
$\Bbb
C^{r_2(\cy)}$ under the canonical  embedding \cite{S,4.2} of $\cy.$ Here
$r_2(\cy)$ is
the number of pairs of complex conjugate embeddings of $\cy$ in $\Bbb C.$ So
$\frac{1}{p^2}
\cy$ also maps to a discrete subgroup.  $Z_p(f,n)$ is represented  by matrices
whose entries
lie in $\frac{1}{p^2} \cy$ and it is an isometry  of the form $<\ ,\  >_\Si$
which is
positive definite under each complex embedding. It follows that each column in
one of these
matrices under each complex embedding must have specified norm with respect to
$<\ ,\
>_\Si$. There can only be a finite number of vectors which meet these
conditions. Thus there
are only a finite number of possible matrices which can describe isometries in
the image of
this representation.
  \qed \enddemo

Funar \cite{F} has proved that the above corollary is false for surfaces of
higher genus.
The reason the above proof does not extend to higher genus or colored points is
that
\cite{BHMV} the form on $<\ ,\  >_\Si,$ is not, in general, positive definite
under all
complex embeddings .

The figure eight knot is a genus one fibered knot.  Thus zero framed surgery
along this knot
is a fiber bundle over a circle with fiber a torus. In \cite{G1,G2} we studied
the map,
$Z_p(F8),$ induced by the monodromy of this bundle for $p \le 20.$  We found
that
in every case  $Z_p(F8)$ had finite order, even though the monodromy itself is
nonperiodic. Here is the list of the periods for $p$ from 3 to 20: 1, 1, 10, 6,
4, 3, 12,
30, 5, 12, 14, 12, 20, 12, 18, 12, 9, 60. We conjectured that $Z_p(F8)$ would
have finite
period for all $p.$ Our attempts to prove it lead to this paper. Of course it
follows from
the Corollary whose truth was unknown to us.

We mention some other related results. Wright has studied the projective
version of these
representations for $p=8$ \cite{Wr1}(with no colored points). One result of her
work is that
the image of this representation is finite. In \cite{Wr2}, Wright studies the
$p=12$ case
for genus two surfaces (without colored points) and can again conclude the
image of this
representation is finite.

In \S 1, we prove Theorem 1. This result
follows easily
from the integrality of an associated quantum invariant of closed manifolds
with links when $p$ is either an
odd prime, or  twice an odd prime. For manifolds without links this is due to
H. Murakami \cite{M1}
\cite{M2}, although
he states the result only for rational homology spheres. Masbaum and Roberts
\cite{MR1} have
given an elegant proof of this result  including the case where the manifolds
contain colored links.

In \S 2, we prove Theorem 2,
 by studying Jeffrey's explicit formula \cite{J} for a related representation
on $(\smallmatrix a & b
\\ c& d\endsmallmatrix)\in SL(2,\bz),$ in terms of $a,$ $b,$ $c,$ and $d.$

\head \S1 Proof of Theorem 1 \endhead

Recall $\eta_p= u_p A_p^3 \frac{A_p^2-A_p^{-2}}{p} G_p(A_p),$ where
$G_p(A_p)=\frac 1 2
\sum_{m=1}^{2p} (-A_p)^{-m^2}\in \bz[A_p]$ \cite{BHMV, MR1}. Recall our
$u_p$ plays the role of $\k_p^3.$
So $\eta_p\in \frac{1}{p}\cy.$
We also note that
$\eta_p^2 =-\frac{(A_p^2-A_p^{-2})^2}{p}\in \frac{1}{p}\cy.$
Thus $\eta_p^g \in { \left(\frac 1 p \right)
}^{[\frac{g+1}{2}]}\cy.$

Let $r$ be an odd prime. Let $p$ be either $r$ or $2r.$
Let $M$ be a closed morphism in $\Cal C_p.$  As in
\cite{MR1}, define $\Cal  I_p(M)=\eta_p^{-1} Z_p(M).$
If $n(M)=0$, then by \cite{MR1} $\Cal  I_p(M)\in \bz[A_p].$ This is stated for
the case that
the colored fat graph is a colored framed link. However it follows immediately
for colored
fat graphs as well since  we may expand a graph into  linear combination of
links over
$\bz[A_p] ,$ as the idempotents  are defined over  $\bz[A_p]$ \cite{MR1}.
Changing $n(M)$ multiplies
$I_p(M)$ by a power of $u_p.$ So if we allow $n(M)$ to be nonzero, we have
$\Cal  I_p(M)\in \cy.$
Thus $<M>_p= \eta_p I_p(M) \in \eta_p \cy.$

\proclaim {Lemma 1} Suppose $\Si$ is  a connected nonempty object in $\Cal C_p$
of genus  $g.$   If $N$ is endomorphism in $\Cal C_p$ from  $\Si$ itself $\Si,$
 the matrix for $Z_p(N)$ with respect
to $\Cal B$ has entries in ${\eta_p}^g \cy\subset {\left(\frac 1
p\right)}^{[\frac{g+1}{2}]}\cy.$ \endproclaim

\demo{Proof}

By \cite{BHMV,4.11},
each $< \b_i,\b_j>_\Si$ is $\delta_i^j \eta_p^{1-g}$ times a unit  in $\cy.$
Perhaps the special case where $\Si$  is a 2-sphere with no colored points
should be discussed
separately. In this case $\Cal G$ is empty,
there is only one coloring $\b,$ and
$<\b,\b>_\Si= \eta_p.$

Let $s_i$ denote $(< \b_i,\b_i>_\Si)^{-1}.$ Thus $s_i\in \eta_p^{g-1}\cy.$
Then $$Z_p(N)(\b_i)=  \sum_{j} s_j < (H,\b_i)\cup N \cup (-H,\b_j)>_p (\b_j).$$
It follows that the matrix for $Z_p(N)$ with respect to $\Cal B$ has entries in
$\eta_p^g\cy.$ \qed \enddemo

If $\Si'$ is another connected object in $\Cal C_p,$ let $\Cal B'$ be a basis
for
$V_p(\Si)$
described as above.
The same proof also proves:

\proclaim {Lemma 2}  For every morphism $N$ from $\Si'$ to $\Si$, the matrix
for $Z_p(N)$ with respect to the bases $\Cal B'$  and $\Cal B$ has entries in
$\eta_p^g \cy.$ \endproclaim

Let $g'$ denote the genus for $\Si'.$ We remark that the roles of the genus of
$g$ and $g'$ are not symmetric. If $\tilde N$ from $\Si$ to $\Si'$ is the
morphism
obtained by reversing the orientation on $N$ and the choice of incoming and
outgoing manifold, then
$Z_p(\tilde N)$ will, with respect to above bases,  have
entries in $\eta_p^{g'} \cy,$ but $Z_p(N)$ may not.

\head \S2 Proof of Theorem 2  \endhead

Assume $p=2r.$   Let $\Si$ be the  boundary of a fixed solid handlebody $H$ of
genus one in $S^3.$ We assign
it the Lagrangian subspace spanned by the (homology class of) a meridian.  We
may let $\Cal
G$ be a core of the solid torus fattened with framing zero.  For $0\le l\le
r-2$ let $b_l$
denote $H$ with $\Cal G$ colored $l.$ Let $\b_l$ denote the element  of
$V_p(\Si)$
represented by $b_l.$ Then $\Cal B =\{\b_l\}_{0\le l\le r-2}.$

We identify the mapping class group of a torus without colored points with
$SL(2,\bz)$ by
looking at the map on homology with respect to a meridian, longitude basis. Let
$S=(\smallmatrix
0&-1\\1 &0\endsmallmatrix),$ $T=(\smallmatrix 1&1\\0 &1\endsmallmatrix).$

{}From the framing relation, one  has
 $\rho_p(T,0)(\b_l)= \mu_l \b_l$ where $\mu_l=(-A_p)^{l^2+2l}=
(-A_p)^{(l+1)^2-1}.$  Also
$<\rho_p(C(S,0))(\b_l),\b_j>_{\Si}$ is $< \ >_p$ of $S^3$ weighted zero
containing a zero framed Hopf link with
components colored $l$ and $j,$ which we denote by $h_{l,j}.$ Since $\Cal B$ is
orthonormal, we have $\rho_p(C(S,0))$ with respect to the basis $\Cal B$   is
given by the matrix $h_{l,j}.$ According to Morton and Strickland \cite{MS},
the bracket
evaluation of this Hopf link is $(-1)^{l+j}[(l+1)(j+1)].$  Here $[n]$ denotes
$\frac{ A_p^{
2n} -A_p^{ -2n} } { A_p^{ 2} -A_p^{ -2} } .$ We remark that the derivation in
\cite{MS} may be
mimiced completely on the skein theory level without reference to
representation theory.
Thus  $h_{l,j}= \eta_p (-1)^{l+j}[(l+1)(j+1)].$
Also one has that $\rho_p(f,n)=(u_p)^n \rho_p(f,0).$

It is now convenient to consider the basis $\Bbb B$ with elements: $\b'_l=
(-1)^{l-1}
b_{l-1}$ for $1\le l\le r-1.$ With respect to $\Bbb B,$ $\rho_p(T,0)$ is given
by $ \hat
S=\d_l^j (-A_p)^{l^2-1},$ and $\rho_p(S,0)$ is given by $\hat S =\eta_p [l j].$

We  adopt the notation
$\z_{t}=e^{\frac{2 \pi i}{t}}.$ Consider the embedding \cite{MR2, note p.134}
$\psi:\bq[\l_p] \rightarrow \bc$ which sends  $A_p$ to $-\a$ where $\a$ denotes
$\z_{2p}$ and  sends $u_p$ to $\z_8^3\a^{-3}.$
Then $\psi (\eta_p)= - \frac{a^2-a^{-2}}{\sqrt{2r}} i.$ Also $\psi (\hat T) =
\d_l^j \a^{l^2-1},$
and $\psi (\hat S) = \frac {-i}{ \sqrt{2r} }  \left( \a^{lj}-\a^{-lj} \right).$
Let $s=4r$ if $r$ is even and $s=8r$ if $r$ is odd.  Then
$\bz[\z_{s}]=\psi(\bz[\l_{p}]).$

Jeffrey \cite{J} studies a matrix representation $\Cal R_p$ of $SL(2,\bz),$
coming from
conformal field theory, which Witten uses in \cite{W}. This is defined by
$\Cal R_p (S)=\psi (\hat S).$ and $ \Cal R_p (T)= \a \z_8^{-1}\psi (\hat T) =
\d_l^j \z_8 \a^{l^2}.$
Since $\a^{-1} \z_8 \in \bz[\z_{s}],$  the theorem will be proved if we can
show that
$p \Cal R_p(\smallmatrix a & b \\ c&
d\endsmallmatrix)\prec \bz[\z_{s}].$  Here we let the symbol $\prec$ mean ``
has all
entries lying in.''
Jeffery's
expression  for $\Cal R (\smallmatrix a & b \\ c& d\endsmallmatrix)$ \cite{J,
2.7} is

$$\Cal R_p(\pmatrix a & b \\ c& d\endpmatrix )_{j,l}=
-\frac{i \z_8^{-\Phi(U)} \text{sign}(c)}{\sqrt{p|c|}}
\z_{4rc}^{dl^2}
\sum\Sb  \g \pmod{2rc}\\ \g=j \pmod{2r} \endSb
\z_{4rc}^{a \g^2}\left( \z_{2rc}^{\g l}- \z_{2rc}^{-\g l}\right)$$

\noindent $\Phi$ denotes the Rademacher $\Phi$-function with values in  $\bz.$
 Thus one
has that

$$\Cal R_p(\pmatrix a & b \\ c& d\endpmatrix)\prec \frac {1}{pc}
\bz[\z_{s|c|}].$$
To see this, we need to observe that $\sqrt{t} \in
\bz[\z_{4t}]$ using a quadratic Gauss sum. We also have:

$$\Cal R_p(\pmatrix c & d \\ -a& -b\endpmatrix)\prec \frac {1}{a \sqrt {p}}
\bz[\z_{s|a|}]\quad
\text{and} \quad
\Cal R_p(\pmatrix 0 & -1 \\ 1& 0\endpmatrix) \prec \frac {1}{\sqrt{p}}
\bz[\z_{s}].$$
As
$(\smallmatrix a & b \\ c& d\endsmallmatrix )= (\smallmatrix 0 & -1 \\ 1&
0\endsmallmatrix) (\smallmatrix c & d \\
-a& -b\endsmallmatrix),$ we  have:

$$\Cal R_p(\pmatrix a & b \\ c& d\endpmatrix)\prec \frac
{1}{p a} \bz[\z_{s|a|}].$$

As $a$,$c$ are relatively prime, $\bq[\z_{s|a|}] \cap
\bq[\z_{s|c|}]=\bq[\z_{s}].$
So $\Cal
R_p(\smallmatrix a & b \\ c& d\endsmallmatrix)\prec \bq[\z_{s}].$ So $p a \Cal
R_p(\smallmatrix a & b \\
c& d\endsmallmatrix)\prec \bz[\z_{s}],$ and $p c \ \Cal R_p(\smallmatrix a & b
\\
c& d\endsmallmatrix)\prec
\bz[\z_{s}].$  Again since $a$,$c$ are relatively prime, $p \Cal
R_p(\smallmatrix a & b \\ c&
d\endsmallmatrix)\prec \bz[\z_{s}].$

  \enddocument \vfill\eject \end